\documentclass[12pt]{article}

\usepackage{fullpage}
\usepackage{hyperref}
\hypersetup{colorlinks,
      linkcolor=blue,
      citecolor=blue,
      urlcolor=blue}
\usepackage{amssymb,amsmath,amsthm}
\usepackage[capitalise]{cleveref}
\usepackage{thm-restate}
\usepackage{xspace}
\usepackage{tikz}
\usepackage{authblk}
\usepackage{subcaption}
\usetikzlibrary{
    decorations.markings,
    decorations.pathmorphing,
    decorations.text,
}

\newtheorem{corollary}{Corollary}
\newtheorem{proposition}{Proposition}
\newtheorem{conjecture}{Conjecture}

\newtheorem{problem}{Problem}

\newcommand{\N}{{\mathbb N}}
\newcommand{\Z}{{\mathbb Z}}
\newcommand{\Real}{{\mathbb R}}
\newcommand{\tw}{{\sf tw}\xspace}
\newcommand{\Ocal}{{\mathcal O}\xspace}
\newcommand{\NP}{{\sf NP}\xspace}
\DeclareMathOperator{\Mad}{Mad}
\DeclareMathOperator{\BBC}{BBC}
\DeclareMathOperator{\CBC}{CBC}
\DeclareMathOperator{\Av}{Ad}
\renewcommand{\epsilon}{\varepsilon}
\renewcommand{\phi}{\varphi}
\renewcommand{\emptyset}{\varnothing}
\let\le\leqslant
\let\ge\geqslant
\let\leq\leqslant
\let\geq\geqslant

\title{Backbone colouring of chordal graphs\footnote{This work has been founded by JST as part of ASPIRE, Grant Number JPMJAP2302, by CNPq (Brazil) projects 313153/2021-3 and 404613/2023-3, by Inria Associated Team CANOE and project CAPES-Cofecub 49587PE (Brazil-France) and ANR-17-EURE-0004 ``EUR DS4H Investments in the Future''.}}

\author[1]{J\'ulio Ara\'ujo}
\author[2]{Nicolas Nisse}
\author[3]{Lucas Picasarri-Arrieta}

\affil[1]{Departamento de Matem\'atica,  Universidade Federal do Cear\'a, Fortaleza, Brazil}
\affil[2]{{Université Côte d'Azur, CNRS, Inria, I3S, Sophia Antipolis, France}}
\affil[3]{National Institute of Informatics, Tokyo, Japan}

\date{}
\begin{document}

\maketitle
\begin{abstract}
A proper $k$-colouring of a graph $G=(V,E)$ is a function $c: V(G)\to \{1,\ldots,k\}$ such that $c(u)\neq c(v)$ for every edge $uv\in E(G)$. The chromatic number $\chi(G)$ is the minimum $k$ such that there exists a proper $k$-colouring of $G$.
Given a spanning subgraph $H$ of $G$, a $q$-backbone $k$-colouring of $(G,H)$ is a proper $k$-colouring $c$ of $G$ such that $\lvert c(u)-c(v)\rvert \ge q$ for every edge $uv\in E(H)$. The $q$-backbone chromatic number $\BBC_q(G,H)$ is the smallest $k$ for which there exists a $q$-backbone $k$-colouring of $(G,H)$. In their seminal paper, Broersma et al.~\cite{BFGW07} ask whether, for any chordal graph $G$ and any spanning forest $H$ of $G$, we have that $\BBC_2(G,H)\leq \chi(G)+\Ocal(1)$.

In this work, we first show that this is true as long as $H$ is bipartite and $G$ is an interval graph in which each vertex belongs to at most two maximal cliques. We then show that this does not extend to bipartite graphs as backbone by exhibiting a family of chordal graphs $G$ with spanning bipartite subgraphs $H$ satisfying $\BBC_2(G,H)\geq \frac{5\chi(G)}{3}$.
Then, we show that if $G$ is chordal and $H$ has bounded maximum average degree (in particular, if $H$ is a forest), then $\BBC_2(G,H)\leq \chi(G)+\Ocal(\sqrt{\chi(G)})$. We finally show that $\BBC_2(G,H)\leq \frac{3}{2}\chi(G)+\Ocal(1)$ holds whenever $G$ is chordal and $H$ is $C_4$-free.
\end{abstract}


\section{Introduction}
\label{sec:intro}
Let $G = (V,E)$ be a graph. Given a positive integer $k$, we denote the set $\{n\in \N\mid 1\leq n\leq k\}$ by $[k]$. A \emph{proper $k$-colouring of $G$} is a function $c:V(G)\rightarrow [k]$ such that $c(u)\neq c(v)$ holds for every edge $uv\in E(G)$. We say $G$ is \emph{$k$-colourable} if there exists a proper $k$-colouring of $G$. The \emph{chromatic number} of $G$, denoted by $\chi(G)$, is the smallest $k$ for which $G$ is $k$-colourable. We say $G$ is \emph{$k$-chromatic} if $\chi(G) = k$. The {\sc Vertex Colouring Problem} consists of determining $\chi(G)$, for a given graph $G$. It is a well-known $\NP$-hard problem~\cite{Karp72} and one of the most studied problems on Graph Theory~\cite{JT95, MR01}.

The {\sc Vertex Colouring Problem} models several practical applications, frequency assignment problems being perhaps the most famous ones~\cite{AHK+07}. There are several variations of the {\sc Vertex colouring Problem} that were defined in order to model the specific constraints of the practical applications related to frequency assignment in networks. Broersma et al.~\cite{BFGW07,Broersma2003} defined the {\sc Backbone Colouring Problem} to model the situation where certain channels of communication are more demanding than others.

Formally, given a graph $G$, a spanning subgraph $H$ of $G$, called the \emph{backbone} of $G$, and two positive integers $q$ and $k$, a \emph{$q$-backbone $k$-colouring of $(G,H)$} is a proper $k$-colouring $c$ of $G$ for which $\lvert c(u)-c(v)\rvert \ge q$ holds for every $uv\in E(H)$. The \emph{$q$-backbone chromatic number of $(G,H)$}, denoted by $\BBC_q(G,H)$, is the minimum $k$ for which there exists a $q$-backbone $k$-colouring of $(G,H)$. The {\sc Backbone Colouring Problem} consists of determining $\BBC_q(G,H)$~\cite{BFGW07,Broersma2003}. 

Observe that, if $H$ is edgeless, then $\BBC_q(G,H) = \chi(G)$; hence computing $\BBC_q(G,H)$ is an $\NP$-hard problem. It has also been proved that the same holds even if $G$ belongs to particular graph classes such as planar graphs and $H$ is a Hamiltonian path or a matching  of $G$~\cite{HKLT14}.

Note that if $f$ is a proper $k$-colouring of $G$, then the function $g: V\to [q\cdot k-q+1]$ defined by $g(v)=q\cdot f(v)-q+1$ is a
$q$-backbone colouring of
$(G,H)$, for any spanning subgraph $H$ of $G$. Moreover it is well-known that
 if $G=H$ and $f$ is a proper $\chi(G)$-colouring of $G$, this $q$-backbone colouring $g$ of
$(G,H)$ is optimal. Therefore, since $\BBC_q(H,H)\leq \BBC_q(G,H)$ and $\BBC_q(G,H)\leq\BBC_q(G,G)$, we have
\begin{equation}\label{eq:trivialbounds}
q\cdot \chi(H)-q+1 \leq \BBC_q(G,H)\leq q\cdot \chi(G)-q+1.  
\end{equation}

In this work, we give a collection of results on the $q$-backbone chromatic number of a given pair $(G,H)$. However, in order to present results in the literature, let us briefly mention other variants of Backbone Colourings that are found in the literature.

\paragraph{Related work.}  A natural and well-studied variant is the one that imposes a circular metric on the colours.  We can see $\Z_k$\footnote{Whenever we refer to $\Z_k$, we mean the group of integers modulo $k$, also denoted by $\Z/k\Z$.} as
a cycle of length $k$ with vertex set $\{1, \dots, k\}$ together with the graphical distance
$|\cdot |_k$. Then $|a-b|_k\geq q$ if and only if $q\le |a-b|\le
k-q$.
A {\it circular $q$-backbone $k$-colouring} of $(G,H)$ is a mapping $f : V(G)\to  \Z_k$ such that $c(v)\neq c(u)$, for each edge $uv\in E(G)$, and $q\le |c(u)-c(v)|\le
k-q$ for each edge $uv\in E(H)$.
The {\it circular $q$-backbone chromatic number} of a graph pair $(G,H)$, denoted $\CBC_q(G,H)$, is the minimum $k$ such that $(G,H)$ admits a circular $q$-backbone $k$-colouring.

Note that if $f$ is a circular $q$-backbone $k$-colouring of $(G,H)$, then $f$ is also a $q$-backbone $k$-colouring of $(G,H)$. On the other hand, observe that 
a $q$-backbone $k$-colouring $f$ of $(G,H)$ is a  circular $q$-backbone $(k+q-1)$-colouring of $(G,H)$.
Hence for every graph pair $(G,H)$, where $H$ is a spanning subgraph of $G$, we have
\begin{equation}
\BBC_q(G,H) \leq \CBC_q(G,H) \leq \BBC_q(G,H)+q-1.\label{eq-triv}
\end{equation}
Combining Inequalities~(\ref{eq:trivialbounds}) and~(\ref{eq-triv}), we observe that 
\begin{equation}
q\cdot \chi(H)-q+1\leq \CBC_q(G,H) \leq q\cdot\chi (G).\label{eq-2chi}
\end{equation}

One can also find a list variant of Backbone Colourings~\cite{HK.14} in the literature and, more recently, Araujo et al.~\cite{ACT24} introduced a variant for oriented backbones. 

\medskip

Let us consider the upper bound provided by~\eqref{eq:trivialbounds} when $q=2$. If $G$ is planar, by the Four Colour Theorem~\cite{appel1977,appel1977a} we have that $\BBC_2(G,H)\leq 7$. If $G$ is a chordal graph, since it is perfect, we know that $\BBC_2(G,H)\leq 2\chi(G)-1 = 2\omega(G)-1$.
In their seminal work, Broersma et al.~\cite{BFGW07} raised two questions that, to the best of our knowledge, are still open.

\begin{problem}[\!\!~\cite{BFGW07}]\label{prob:planar-forest}
If $G$ is planar and $H$ is a spanning forest of $G$, then is it true that $\BBC_2(G,H)\leq 6$?
\end{problem}

\begin{problem}[\!\!~\cite{BFGW07}]\label{prob:chordal-forest}
Does there exist a positive integer $c$ such that, for any chordal graph $G$ and any spanning forest $H$ of $G$, we have that $\BBC_2(G,H)\leq \omega(G)+c$?    
\end{problem}

Concerning \Cref{prob:planar-forest}, the authors also asked for a proof that $\BBC_2(G,H)\leq 7$, given a planar graph $G$ and a spanning forest $H$ of $G$, without using the Four Colour Theorem.
Many authors in the literature considered that the answer to \Cref{prob:planar-forest} should be affirmative and it is treated as a conjecture by most of the related works. This conjecture would be tight even if $H$ is a Hamiltonian path, as there are examples of a planar graph $G$ and Hamiltonian path $H$ in $G$ for which $\BBC_2(G,H)= 6$~\cite{BFY.03}. \Cref{prob:planar-forest} has also a version proposed by the same authors when $H$ is restricted to be a matching, in which case the upper bound is conjectured to be 5~\cite{BFY.03}. There is also another variation of \Cref{prob:planar-forest} for circular backbone colourings where the upper bound is 7~\cite{HKLT14}.

Some effort to solve \Cref{prob:planar-forest} and its variations has been done by different authors in the last decade. Campos et al.~\cite{CHSS13} confirmed that \Cref{prob:planar-forest} has a positive answer if $H$ is a tree of diameter at most 4. For larger values of $q$, Havet et al.~\cite{HKLT14} proved that, if $G$ is planar and $H$ is a forest in $G$, then $\BBC_q(G,H)\leq q+6$ and that this is tight for $q=4$. Regarding circular backbone colouring, Havet et al.~\cite{HKLT14} proved that if $G$ is planar and $H$ is a forest in $G$, then $\CBC_q(G,H)\leq 2q+4$ and conjectured that this bound is not tight. This was latter verified by~\cite{AAC+22}, that also presented partial results for the circular variant of \Cref{prob:planar-forest}. Araujo et al.~\cite{AHS.15} proved bounds when $G$ is a planar graph without cycles of length 4 or 5 and $H$ is a spanning forest of $G$ in the circular variant of \Cref{prob:planar-forest}.

On the other hand, \Cref{prob:chordal-forest} was not much studied in the literature. In case $G$ is chordal and $H$ is a Hamiltonian path of $G$, Broersma et al.~\cite{BFGW07} observed that $\BBC_2(G,H)\leq \chi(G) +4$. Actually, the technique they used can be easily extended to prove that if $G$ is chordal and $H$ has degree at most $c$, then $\BBC_2(G,H)\leq \omega(G)+2c$.
A particular case of chordal graphs that has received some attention is the one of split graphs~\cite{salman2006lambda,broersma2009backbone,turowski2015}. Salman~\cite{salman2006lambda} presented tight upper bounds for $\BBC_q(G,H)$ in function of $\chi(G)$ and $q$, when $G$ is split and $H$ is a forest of stars. The authors in~\cite{broersma2009backbone} proved similar bounds when $G$ is split and $H$ is a matching in $G$. They also presented complexity results for the computation of $\BBC_q(G,H)$ for arbitrary $G$ and $H$ being a matching. Turowski~\cite{turowski2015} showed bounds when $G$ is complete and $H$ belongs to some classes, and also presented a polynomial-time algorithm to compute $\BBC_2(G,H)$ when $G$ is split and $H$ is a matching in $G$.

\paragraph{Our contribution.} In this work, we present some partial results towards the resolution of \Cref{prob:chordal-forest} by possibly constraining $G$ to some subclasses of chordal graphs and also releasing the constraint of $H$ to be a forest.
We first consider the case of $H$ being bipartite. When $G$ is an interval graph in which each vertex belongs to at most two maximal cliques, we give a positive answer to~\Cref{prob:chordal-forest}.

\begin{restatable}{theorem}{intervaltwocliques}\label{thm:interval2cliques}
Let $G$ be an interval graph such that every vertex belongs to at most $2$ maximal cliques and let $H$ be a bipartite spanning subgraph of $G$. Then, $\BBC_2(G,H)\leq \omega(G)+3$.
\end{restatable}

We then show that Theorem~\ref{thm:interval2cliques} does not extend to the more general case of $G$ being chordal by exhibiting a family of chordal graphs $G$ with spanning bipartite subgraphs $H$ satisfying $\BBC_2(G,H)\geq \frac{5}{3}\omega(G)$.
\begin{restatable}{proposition}{prochordalbip}
    \label{prop:example_chordal_bip_53w}
    For infinitely many values of $\omega$, there exists a chordal graph $G$ and a spanning bipartite subgraph $H$ of $G$ such that $\omega(G) = \omega$ and $\BBC_2(G,H)\geq \frac{5}{3}\omega$.
\end{restatable}

We then consider the case of $H$ having bounded maximum average degree.
Recall that $\Av(H') = \frac{2|E(H')|}{|V(H')|}$ is the {\it average degree} of $H$ and that is $\Mad(H) = \max \{\Av(H')\mid H'\subseteq H\}$ is its \emph{maximum average degree}. Note also that forests have maximum average degree smaller than $2$.
When $\Mad(H)\leq d$, we prove that $\displaystyle \BBC_2(G,H)\leq (1+o(1))\cdot \omega(G)$.

\begin{restatable}{theorem}{thmmad} \label{thm:mad}
    Let $d\in\Real^*_+$. Let $G$ be any chordal graph and let $H$ be any subgraph of $G$ with $\Mad(H)\leq d$. Then, $\BBC_2(G,H)\leq \omega(G)+2 \sqrt{d \cdot \omega(G)} +3d$.
\end{restatable}

We finally show that $\BBC_2(G,H)\leq \frac{3}{2}\omega(G)+\Ocal(1)$ holds whenever $G$ is chordal and $H$ is $C_4$-free, which is thus halfway between~\Cref{prob:chordal-forest} and the trivial bound $\BBC_2(G,H)\leq 2\chi(G)-1$. Recall that $H$ is $C_4$-free if $H$ does not contain $C_4$ as a subgraph (not necessarily induced).

\begin{restatable}{theorem}{thmcfourfree}
\label{thm:3half}
Let $G$ be a chordal graph and $H$ be a $C_4$-free spanning subgraph of $G$. Then, $\BBC_2(G,H)\leq \frac{3}{2}\omega(G)+4$.
\end{restatable}

A \textit{tree-decomposition} of $G$ is a pair $(T,\mathcal{X})$ where $T=(I,F)$ is a tree, and $\mathcal{X}=(B_i)_{i\in I}$ is a family of subsets of $V(G)$, called \textit{bags} and indexed by the vertices of $T$, such that:
\begin{enumerate}
    \item each vertex $v\in V$ appears in at least one bag, \textit{i.e.} $\bigcup_{i\in I} B_i= V$,
    \item for each edge $e = xy \in E$, there is an $i\in I$ such that $x,y \in B_i$, and 
    \item for each $v\in V$, the set of nodes indexed by $\{ i \mid i\in I, v\in B_i\}$ forms a subtree of $T$.
\end{enumerate}

The \textit{width} of a tree decomposition is defined as $\max_{i\in I} \{|B_i| -1\}$. The \textit{treewidth} of $G$, denoted by $\tw (G)$, is the minimum width of a tree-decomposition of $G$.
It is well-known that every graph $G$ is a subgraph of a chordal graph $G'$ with $\omega(G') = \tw(G) + 1$ (see for instance~\cite[Corollary~12.3.12]{diestel2017}). Hence the following is a direct consequence of Theorems~\ref{thm:mad} and~\ref{thm:3half}.
\begin{corollary}
    \label{cor-tw}
    Let $G$ be a graph and $H$ be a subgraph of $G$, then:
    \begin{enumerate}
        \item $\BBC_2(G,H) \leq \tw(G) + \Ocal(\sqrt{\Mad(H) \cdot \tw(G)})$, and
        \item $\BBC_2(G,H) \leq\frac{3}{2}\tw(G)+\frac{11}{2}$ if $H$ is $C_4$-free.
    \end{enumerate}
\end{corollary}

\section{Our Contributions}
\label{sec:contributions}

For basic notions and terminology on Graph Theory not explicitly defined in this section, the reader is referred to~\cite{BM.book}.  
Let $G$ be a graph and $X\subseteq V(G)$ be any subset of vertices. We denote by $N_G[X]$ the union of closed neighbourhoods in $X$, that is $N_G[X] = \bigcup_{x\in X} N[x]$, where $N[x]$ denotes $\{x\}\cup N(x)$.

A {\it chordal graph} is a graph without any induced cycle of length at least $4$. A graph $G=(V,E)$ is an {\it interval graph} if there exists a set of intervals $\mathcal{I}$ on the real line and an injection $\iota\colon V\to \mathcal{I}$ such that, for every pair of vertices $u,v\in V$, $uv\in E$ if and only if $\iota(u)\cap \iota(v) \neq \emptyset$. It is well known that interval graphs are chordal (see~\cite[Exercise 5.42]{diestel2017}), and that the converse is not true.

A \textit{path-decomposition} is a tree-decomposition  $(P,\mathcal{X})$ with the extra property that $P$ is a path. A clique $C$ of a graph $G$ is {\it maximal} if every vertex of $V(G)\setminus C$ is non-adjacent to at least one vertex in $C$.
We skip the proof of the following well-known property of interval graphs (see for instance~\cite[Exercise~12.29]{diestel2017}).

\begin{proposition}
    \label{prop:path_dec_interval}
    Every interval graph $G$ admits a path-decomposition $(P,\mathcal{X}=(X_1,\dots,X_\ell))$ such that, for every $i\in [\ell]$, $X_i$ is a maximal clique of $G$.
\end{proposition}

\subsection{Bipartite backbones}
The goal of this section is to prove \Cref{thm:interval2cliques} and \Cref{prop:example_chordal_bip_53w}.

\intervaltwocliques*

\begin{proof}
In this proof, since we will consider circular permutations of colours modulo $\omega(G)$, let us assume for better readability that the colours belong to $\{0,\cdots,\omega(G)-1\}$ and that all arithmetic operations with colours must be understood modulo $\omega(G)$. A {\it circular interval}, denoted by $[a,b]$, is the interval $[a,a+1,\cdots,b-1,b]$ if $0 \leq a \leq b <\omega(G)$ and it is equal to $[a,a+1,\cdots,\omega(G)-1,0,\cdots,b]$ if $0 \leq b < a <\omega(G)$. With a slight abuse of notation, when it is convenient, we consider circular intervals as sets, otherwise as sequences of consecutive non-negative integers (modulo $\omega(G)$). In particular, we consider that circular intervals are always represented in clockwise order.

Let $(A,B)$ be a bipartition of $H$, and $(X_1,\cdots,X_{\ell})$ be a path-decomposition of $G$ such that $X_i$ is a maximal clique for every $1 \leq i \leq \ell$, the existence of which is guaranteed by Proposition~\ref{prop:path_dec_interval}. Without loss of generality, we assume that each $X_i$ has size exactly $\omega(G)$. Indeed, if this is not the case, for each $i\in[t]$ we add $\omega(G)-|X_i|$ new pairwise adjacent vertices dominating $X_i$ and these new vertices belong to $A$ and are isolated vertices in $H$.
In what follows, we further let $X_0 = X_{\ell+1} = \emptyset$.

By assumption, $X_{i-1} \cap X_{i+1}=\emptyset$ for every $1<i<\ell$. Moreover, for every $1 \leq i \leq \ell$, let $A_i=X_{i-1} \cap X_{i}\cap A$, $B_i=X_{i-1} \cap X_{i}\cap B$, $A'_i = X_i \cap A \setminus (X_{i-1} \cup X_{i+1})$ (i.e., $A'_i$ is the set of the simplicial vertices of $X_i \cap A$) and $B'_i = X_i \cap B \setminus (X_{i-1} \cup X_{i+1})$. Note that $A_1=B_1=\emptyset$. Note also that $\{A_i, A'_i, A_{i+1}\}$ is a partition of $X_i \cap A$, possibly with empty parts, while the same happens for $\{B_i, B'_i, B_{i+1}\}$ with respect to $X_i \cap B$, for every $1 \leq i \leq \ell$. See \Cref{fig:notationAiBi}.

\begin{figure}[!ht]
    \centering
    \begin{tikzpicture}
    \newcommand{\rad}{2} 
    \newcommand{\dist}{3pt} 
    \begin{scope}
        \draw (0,0) circle [radius=\rad];
        \draw (-3,0) circle [radius=\rad];
        \draw (3,0) circle [radius=\rad];
        \draw (-6.5,0) edge [dotted, color=red] (6.5,0);
        \draw (-6,-2) arc (-90:90:\rad);
        \draw (6,2) arc (90:270:\rad);

        \draw (-3,1.25) node [color=red] {$A_{i-1}'$};
        \draw (-1.5,0.5) node [color=red] {$A_i$};
        \draw (0,1.25) node [color=red] {$A_i'$};
        \draw (1.5,0.5) node [color=red] {$A_{i+1}$};
        \draw (3,1.25) node [color=red] {$A_{i+1}'$};

        \draw (-3,-1.25) node [color=red] {$B_{i-1}'$};
        \draw (-1.5,-0.5) node [color=red] {$B_i$};
        \draw (0,-1.25) node [color=red] {$B_i'$};
        \draw (1.5,-0.5) node [color=red] {$B_{i+1}$};
        \draw (3,-1.25) node [color=red] {$B_{i+1}'$};
    \end{scope}
    \end{tikzpicture}        
    \caption{Representation of sets used in the proof of \Cref{thm:interval2cliques}}
    \label{fig:notationAiBi}
\end{figure}
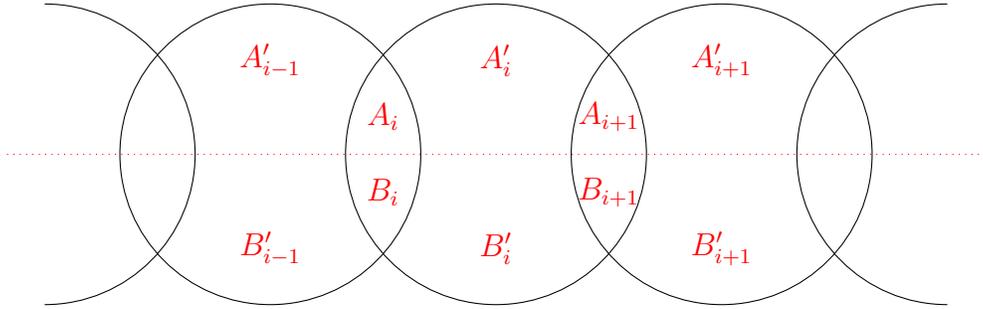

We start by building a proper $\omega(G)$-colouring $\phi:V(G) \rightarrow \{0,\cdots,\omega(G)-1\}$ of $G$ with the following properties: for every $1 \leq i \leq \ell$,
\begin{itemize}
    \item $\phi(X_i)=\{0,\cdots,\omega(G)-1\}$;
    \item $\phi(A \cap X_i)$ and $\phi(B\cap X_i)$ are two circular intervals that partition $\{0,\cdots,\omega(G)-1\}$.
    \item $\phi(A \cup X_i)$ is the concatenation of the three circular intervals $\phi(A_i),\phi(A'_i)$ and $\phi(A_{i+1})$ (in this order if $i$ is odd and in the reversed order if $i$ is even), and
    \item $\phi(B \cup X_i)$ is the concatenation of the three circular intervals $\phi(B_i),\phi(B'_i)$ and $\phi(B_{i+1})$ (in this order if $i$ is even and in the reversed order if $i$ is odd).
\end{itemize}

To show that such a colouring $\phi$ exists, we inductively build its restriction $\phi_i$ to $G[\bigcup_{j\leq i}X_j]$ for each $i$ from $1$ to $\ell$ (obviously $\phi_{\ell}=\phi$).

Let $\phi_1$ be defined as follows:
\begin{itemize}
    \item If $X_i \cap A =\emptyset$, then $\phi_1$ colours the vertices of $B'_1$ with colours $0$ to $|B'_1|-1$, and then vertices of $B_2$ with colours $|B'_1|$ to $|B'_1|+|B_2|-1=|X_1|-1=\omega(G)-1$.
    \item Otherwise, $\phi_1$ colours the vertices of $A'_1$ with colours $0$ to $|A'_1|-1$, and then vertices of $A_2$ with colours $|A'_1|$ to $|A'_1|+|A_2|-1=|A \cap X_1|-1$. Then, $\phi_1$ colours the vertices of $B_2$ from colour $|A \cap X_1|$ to $|A \cap X_1|+|B_2|-1$ and then vertices of $B'_1$ from $|A \cap X_1|+|B_2|$ to $|A \cap X_1|+|B_2|+|B'_1|-1=|X_1|-1=\omega(G)-1$.

    In this case, let denote by $v_1$ (resp., by $v'_1$) the vertex of $A \cap X_1$ coloured with colour $0$ (resp., the vertex of $A \cap X_1$ coloured with colour $|A \cap X_1|-1$). Possibly $v_1=v'_1$.
\end{itemize}

Now, let us assume that $\phi_{i-1}$ has been defined recursively for some $1<i\leq \ell$, and that, by induction, there exist colours $x_i,y_i,c_i \in \{0,\cdots,\omega(G)-1\}$ and $c_i \in [x_i,y_i]$, such that:
\begin{itemize}
    \item either $A_i=\emptyset$ and $\phi_{i-1}(B_i)=[x_i,y_i]$, or $B_i=\emptyset$ and $\phi_{i-1}(A_i)=[x_i,y_i]$ (note that $A_i \cup B_i \neq \emptyset$ since $G$ is connected), or 
        \item $\phi_{i-1}(A_i \cup B_i)=[x_i,y_i]$ and  $\phi_{i-1}(A_i)=[x_i,c_i]$ and $\phi_{i-1}(B_i)=[c_i+1,y_i]$, or vice-versa:
            \item $\phi_{i-1}(A_i \cup B_i)=[x_i,y_i]$ and  $\phi_{i-1}(B_i)=[x_i,c_i-1]$ and $\phi_{i-1}(A_i)=[c_i,y_i]$. 
\end{itemize}

We extend $\phi_{i-1}$ to $\phi_i$ as follows. Obviously, for every $v \in \bigcup_{0\leq j< i}X_j$, $\phi_i(v)=\phi_{i-1}(v)$. For the vertices of $X_i \setminus X_{i-1}$, there are the following cases to be considered:
\begin{itemize}
    \item First let us consider the case when $A \cap (X_i \setminus X_{i-1}) = A'_i \cup A_{i+1}=\emptyset$. 
    \begin{itemize}
\item     If $\phi_{i-1}(A_i \cup B_i)=[x_i,y_i]$ is the concatenation of the interval $\phi_{i-1}(A_i)$ and $\phi_{i-1}(B_i)$ (in this order) and possibly $A_i$ or $B_i$ being empty, then $\phi_i$ colours the vertices of $B'_i$ with colours from $y_i+1$ to $y_i+|B'_i|$ and then $\phi_i$ colours the vertices of $B_{i+1}$ with colours from $y_i+|B'_i|+1$ to $y_i+|B'_i|+|B_{i+1}|= x_i-1$ (modulo $\omega(G)$).

(the last equality comes from the fact that $|A_i \cup B_i| = |X_i \cap X_{i-1}| = |[x_i,y_i]|$ and $\omega(G)=|X_i|= |X_i \cap X_{i-1}| + |B'_i|+|B_{i+1}|$.)
\item     Otherwise, $\phi_{i-1}(A_i \cup B_i)=[x_i,y_i]$ is the concatenation of the interval $\phi_{i-1}(B_i)$ and $\phi_{i-1}(A_i)$ (in this order). Then $\phi_i$ colours the vertices of $B_{i+1}$ with colours from $y_i+1$ to $y_i+|B_{i+1}|$ and then $\phi_i$ colours the vertices of $B'_i$ with colours from $y_i+|B_{i+1}|+1$ to $y_i+|B_{i+1}|+|B'_i|=x_i-1$ (modulo $\omega(G)$).
\end{itemize}
\item Otherwise, $A \cap (X_i \setminus X_{i-1}) = A'_i \cup A_{i+1}\neq \emptyset$. 
\begin{itemize}
    \item If $\phi_{i-1}(A_i \cup B_i)=[x_i,y_i]$ is the concatenation of the interval $\phi_{i-1}(A_i)$ and $\phi_{i-1}(B_i)$ (in this order) and possibly $A_i$ or $B_i$ is empty, then $\phi_i$ colours the vertices of $B'_i$ with colours from $y_i+1$ to $y_i+|B'_i|$ and then $\phi_i$ colours the vertices of $B_{i+1}$ with colours from $y_i+|B'_i|+1$ to $y_i+|B'_i|+|B_{i+1}|$. Then, $\phi_i$ colours the vertices of $A_{i+1}$ with colours from $y_i+|B'_i|+|B_{i+1}|+1$ to $y_i+|B'_i|+|B_{i+1}|+|A_{i+1}|$ and then $\phi_i$ colours $A'_i$ with colours from $y_i+|B'_i|+|B_{i+1}|+|A_{i+1}|+1$ to $y_i+|B'_i|+|B_{i+1}|+|A_{i+1}|+|A'_i|=x_i-1 \mod \omega(G)$. See Figure~\ref{fig:coloursforintervals}(b).

    (the last equality comes from the fact that $|A_i \cup B_i| = |X_i \cap X_{i-1}| = |[x_i,y_i]|$ and $\omega(G)=|X_i|= |X_i \cap X_{i-1}| + |X_i \setminus X_{i-1}|$ and $X_i \setminus X_{i-1}=A'_i \cup A_{i+1}\cup B'_i\cup B_{i+1}$.)

    \item If $\phi_{i-1}(A_i \cup B_i)=[x_i,y_i]$ is the concatenation of the interval $\phi_{i-1}(B_i)$ and $\phi_{i-1}(A_i)$ (in this order), then $\phi_i$ colours the vertices of $A'_i$ with colours from $y_i+1$ to $y_i+|A'_i|$ and then $\phi_i$ colours the vertices of $A_{i+1}$ with colours from $y_i+|A'_i|+1$ to $y_i+|A'_i|+|A_{i+1}|$. Then, $\phi_i$ colours the vertices of $B_{i+1}$ with colours from $y_i+|A'_i|+|A_{i+1}|+1$ to $y_i+|A'_i|+|A_{i+1}|+|B_{i+1}|$ and finally $\phi_i$ colours $B'_i$ with colours from $y_i+|A'_i|+|A_{i+1}|+|B_{i+1}|+1$ to $y_i+|A'_i|+|A_{i+1}|+|B_{i+1}|+|B'_i|=x_i-1$ (modulo $\omega(G)$). See Figure~\ref{fig:coloursforintervals}(a).
\end{itemize}
Defined that way, we have that $\phi_i(A \cap X_i)=[\alpha_i,\beta_i]$. Let $v_i$ (resp., $v'_i$) be the vertex of $A \cap X_i$ coloured with colour $\alpha_i$ (resp., with colour $\beta_i$). 
\end{itemize}

   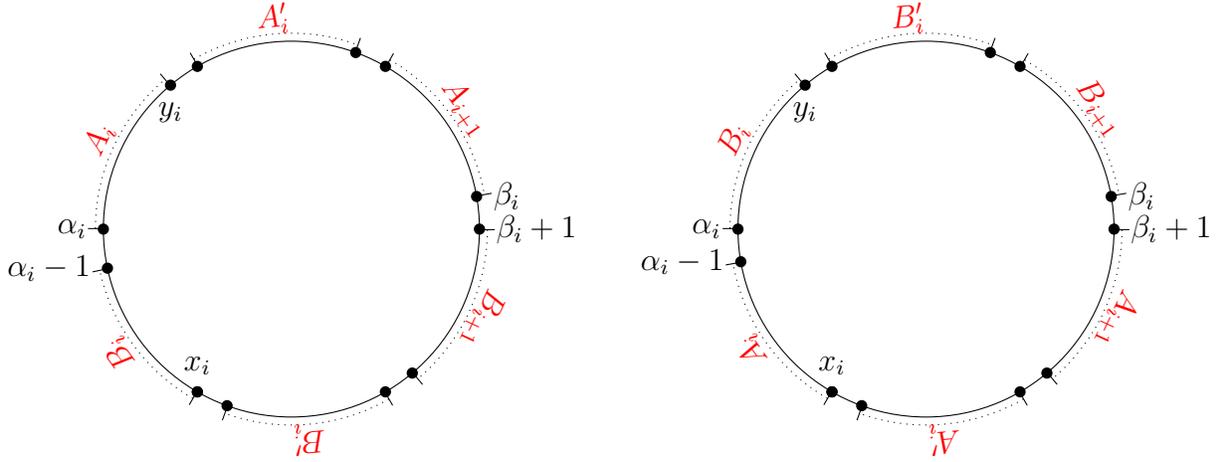
\begin{figure}[!ht]
    \centering
        \begin{subfigure}[b]{.49\textwidth}
            \begin{tikzpicture}
        \newcommand{\rad}{2.5} 
        \newcommand{\dist}{3pt} 
            \begin{scope}
                \clip[draw] circle [radius=\rad];
            \end{scope}
            \draw (192:\rad) node[circle,inner sep=1.5pt,fill, label=left:$\alpha_{i}-1$] {};
            \draw (180:\rad) node[circle,inner sep=1.5pt,fill, label=left:$\alpha_{i}$] {};
            \draw (130:\rad) node[circle,inner sep=1.5pt,fill, label=below:$y_i$] {};
            \draw (120:\rad) node[circle,inner sep=1.5pt,fill] {};
            \draw (70:\rad) node[circle,inner sep=1.5pt,fill] {};
            \draw (60:\rad) node[circle,inner sep=1.5pt,fill] {};
            \draw (10:\rad) node[circle,inner sep=1.5pt,fill, label=right:$\beta_i$] {};
            \draw (0:\rad) node[circle,inner sep=1.5pt,fill, label=right:$\beta_i+1$] {};
            \draw (-50:\rad) node[circle,inner sep=1.5pt,fill] {};
            \draw (-60:\rad) node[circle,inner sep=1.5pt,fill] {};
            \draw (-110:\rad) node[circle,inner sep=1.5pt,fill] {};
            \draw (-120:\rad) node[circle,inner sep=1.5pt,fill] {};
            \draw (-120:\rad) node[circle,inner sep=1.5pt,fill, label=above:$x_i$] {};

            \draw [dotted,|-|, postaction={decorate, decoration={text align={center},raise={1mm},text along path, text={|\color{red}|{$A_{i}$}}}}] (180:\rad cm+\dist) arc (180:130:\rad cm+\dist);
            \draw [dotted,|-|, postaction={decorate, decoration={text align={center},raise={1mm},text along path, text={|\color{red}|{$A_{i}'$}}}}] (120:\rad cm+\dist) arc (120:70:\rad cm+\dist);
            \draw [dotted,|-|, postaction={decorate, decoration={text align={center},raise={1mm},text along path, text={|\color{red}|{$A_{i+1}$}}}}] (60:\rad cm+\dist) arc (60:10:\rad cm+\dist);
            \draw [dotted,|-|, postaction={decorate, decoration={text align={center},raise={1mm},text along path, text={|\color{red}|{$B_{i+1}$}}}}] (0:\rad cm+\dist) arc (0:-50:\rad cm+\dist);
            \draw [dotted,|-|, postaction={decorate, decoration={text align={center},raise={1mm},text along path, text={|\color{red}|{$B_{i}'$}}}}] (-60:\rad cm+\dist) arc (-60:-110:\rad cm+\dist);
            \draw [dotted,|-|, postaction={decorate, decoration={text align={center},raise={1mm},text along path, text={|\color{red}|{$B_{i}$}}}}] (-120:\rad cm+\dist) arc (-120:-168:\rad cm+\dist);
            
            \draw (0,-\rad-.5) node {}; %for caption not to touch labeks
        \end{tikzpicture}
        \caption{Case when $\phi_{i-1}(A_i \cup B_i)=[x_i,y_i]$ is the concatenation of the interval $\phi_{i-1}(B_i)$ and $\phi_{i-1}(A_i)$ (in this order).}
        \label{fig:oddcolouring}
        \end{subfigure}
        \hfill
        \begin{subfigure}[b]{.49\textwidth}
            \begin{tikzpicture}
        \newcommand{\rad}{2.5} 
        \newcommand{\dist}{3pt} 
            \begin{scope}
                \clip[draw] circle [radius=\rad];
            \end{scope}
            \draw (190:\rad) node[circle,inner sep=1.5pt,fill, label=left:$\alpha_i-1$] {};
            \draw (180:\rad) node[circle,inner sep=1.5pt,fill, label=left:$\alpha_{i}$] {};
            \draw (130:\rad) node[circle,inner sep=1.5pt,fill, label=below:$y_i$] {};
            \draw (120:\rad) node[circle,inner sep=1.5pt,fill] {};
            \draw (70:\rad) node[circle,inner sep=1.5pt,fill] {};
            \draw (60:\rad) node[circle,inner sep=1.5pt,fill] {};
            \draw (10:\rad) node[circle,inner sep=1.5pt,fill, label=right:$\beta_i$] {};
            \draw (0:\rad) node[circle,inner sep=1.5pt,fill, label=right:$\beta_i+1$] {};
            \draw (-50:\rad) node[circle,inner sep=1.5pt,fill] {};
            \draw (-60:\rad) node[circle,inner sep=1.5pt,fill] {};
            \draw (-110:\rad) node[circle,inner sep=1.5pt,fill] {};
            \draw (-120:\rad) node[circle,inner sep=1.5pt,fill] {};
            \draw (-120:\rad) node[circle,inner sep=1.5pt,fill, label=above:$x_i$] {};

            \draw [dotted,|-|, postaction={decorate, decoration={text align={center},raise={1mm},text along path, text={|\color{red}|{$B_{i}$}}}}] (180:\rad cm+\dist) arc (180:130:\rad cm+\dist);
            \draw [dotted,|-|, postaction={decorate, decoration={text align={center},raise={1mm},text along path, text={|\color{red}|{$B_{i}'$}}}}] (120:\rad cm+\dist) arc (120:70:\rad cm+\dist);
            \draw [dotted,|-|, postaction={decorate, decoration={text align={center},raise={1mm},text along path, text={|\color{red}|{$B_{i+1}$}}}}] (60:\rad cm+\dist) arc (60:10:\rad cm+\dist);
            \draw [dotted,|-|, postaction={decorate, decoration={text align={center},raise={1mm},text along path, text={|\color{red}|{$A_{i+1}$}}}}] (0:\rad cm+\dist) arc (0:-50:\rad cm+\dist);
            \draw [dotted,|-|, postaction={decorate, decoration={text align={center},raise={1mm},text along path, text={|\color{red}|{$A_{i}'$}}}}] (-60:\rad cm+\dist) arc (-60:-110:\rad cm+\dist);
            \draw [dotted,|-|, postaction={decorate, decoration={text align={center},raise={1mm},text along path, text={|\color{red}|{$A_{i}$}}}}] (-120:\rad cm+\dist) arc (-120:-170:\rad cm+\dist);

            \draw (0,-\rad-.5) node {}; 
        \end{tikzpicture}
        \caption{Case when $\phi_{i-1}(A_i \cup B_i)=[x_i,y_i]$ is the concatenation of the interval $\phi_{i-1}(A_i)$ and $\phi_{i-1}(B_i)$ (in this order).}
        \label{fig:evencolouring}
        \end{subfigure}
        
        \caption{Colours' ordering for vertices in $X_i$ in the proof of \Cref{thm:interval2cliques}}
        \label{fig:coloursforintervals}
    \end{figure}

    The function $\phi$ that we defined above is clearly a $\omega(G)$-proper colouring of $G$ satisfying the announced properties.

To obtain the desired 2-backbone $(\omega(G)+3)$-colouring $\phi'$ of $(G,H)$, let us start from a colouring $\phi'=\phi$. Now, note that the set $X=\bigcup_{1 \leq i \leq \ell} \{v_i,v'_i\} \subseteq A$ induces a set of disjoint paths in $G$. Indeed, each clique $X_i$ contains at most two vertices of $X$ ($v_i$ and/or $v'_i$), and, for $v \in \{v_i,v'_i\}$, if $v \in X_i \cap X_{i-1}$, then $v \in \{v_{i-1},v'_{i-1}\}$. We modify $\phi'$ by recolouring vertices in $X$ in such a way that $\phi'$ induces a proper colouring of $G[X]$ using colours $\{\omega(G)+1,\omega(G)+2\}$.

Note that $\phi'$ now is a proper $(\omega(G)+3)$-colouring of $G$. 
Now, let $\{u,v\} \in E(H)$ and let $1 \leq i \leq \ell$ such that $u \in X_i \cap A$ and $v \in X_i \cap B$. By construction of $\phi$, the colour $\phi(v)=\phi'(v)$ is at distance $2$ from $\phi(w)=\phi'(w)$ for every $w \in X_i \cap A \setminus \{v_i,v_{i+1}\}$. Moreover, $\{\phi'(v_i),\phi'(v_{i+1})\}=\{\omega(G)+1,\omega(G)+2\}$ and $\phi'(v)<\omega(G)$. Hence, the distance between $\phi'(u)$ and $\phi'(v)$ is at least two. So $\phi'$ is a 2-backbone $(\omega(G)+3)$-colouring of $(G,H)$. This concludes the proof.
\end{proof}

In the sequel, we show that the upper bound provided by \Cref{thm:interval2cliques} cannot be extended to the more general case of $G$ being chordal.

\prochordalbip*

\begin{proof}
Let us first build $G$ and $H$ as required. Let $r \in \N$ be any positive integer and define $\omega=3r$. To build $G$, we start with a complete graph $K_{\omega}$. We refer to this initial set of vertices as $K$. Then, for any $X \subseteq K$ of size $|X|=r$, add a complete graph $K_{2r}$ on $2r$ new vertices and make them adjacent to all vertices in $X$. For each $X$, we refer to such $2r$ vertices as $K_X$. It is straightforward to check that $G$ is chordal.

Let $H$ be the bipartite graph induced by the edges linking the vertices in $K$ to the ones in $K_X$, for all $X \subseteq K$ with $|X|=r$. Note that $K$ is an independent set in $H$, as well as the union of the vertices in $K_X$, over all $X\subseteq K$ with $|X|=r$. Thus, $H$ is bipartite. 

Let $\phi$ be a 2-backbone $k$-colouring of $(G,H)$. Note that a vertex $v\in K$ with colour $1<a<k$ forbids not only colour $a$, but also colours $a-1$ and $a+1$ from appearing in its neighbourhood in $H$. Since $K$ induces a complete graph on $G$, in $\phi$ there are $\omega$ distinct colours $1\leq c_1<c_2<\cdots<c_{\omega}\leq k$ in $K$. Let $X$ be the set of vertices $v\in K$ such that $\phi(v) = c_{3p+2}$ for some $p\in \{0,\ldots, r-1\}\}$. Note that $|X| = r$ and each vertex in $X$ forbids a distinct set of exactly three colours from appearing in all vertices in $K_X$ in the 2-backbone $k$-colouring $c$. Consequently, $\phi$ uses at least $3r+2r=\omega+2r$ colours as $K_X$ induces a complete graph on $2r$ vertices. This shows 
that 
$\BBC_2(G,H) \geq \frac{5}{3}\omega$, 
as desired.
\end{proof}

\subsection{Sparse backbones}

This section is devoted to the proof of Theorem~\ref{thm:mad}.

\thmmad*
\begin{proof}
Fix $d\in\Real^*_+$. Let us show a stronger statement that, for every $\epsilon>0$, every chordal graph $G$, and every spanning subgraph $H$ of $G$ with $\Mad(H)\leq d$, we have 
\[ 
\BBC_2(G,H)\leq (1+\epsilon)\omega(G)+c_{\epsilon,d},
\]
where $c_{\epsilon,d}=\max \left\{ \frac{d}{\epsilon}, \frac{d}{2\epsilon}+3d\right\} \leq \frac{d}{\epsilon}+3d$. In particular, by setting $\epsilon=\sqrt{\frac{d}{\omega(G)}}$, we obtain
\begin{align*}
    \BBC_2(G,H)
    \leq \left(1+\epsilon\right)\omega(G) + \frac{d}{\epsilon}+3d
    &\leq \left(1+\sqrt{\frac{d}{\omega(G)}}\right)\omega(G) + \frac{d}{\sqrt{\frac{d}{\omega(G)}}}+3d\\
    &\leq \omega(G)+2 \sqrt{d \cdot \omega(G)} +3d.
\end{align*}

Fix a real $\epsilon>0$ and let $c_{\epsilon,d}$ be as defined above. 
We will show by induction on the order $n=|V(G)|=|V(H)|$ of $G$ that $\BBC_2(G,H) \leq (1+\epsilon)\omega(G) + c_{\epsilon,d}$. 
For better readability, in what follows, we denote $\omega(G)$ by $\omega$ and we fix $k$ to $\lfloor(1+\epsilon)\omega + c_{\epsilon,d}\rfloor$.

Assume first that $n\leq \frac{d}{\epsilon}$. In this case, it is sufficient to prove that $\BBC_2(G,H) \leq 2\omega$ as $2\omega \leq \omega + n \leq k$. We are thus done by~\eqref{eq:trivialbounds}.
Henceforth we assume that $n > \frac{d}{\epsilon}$ and that the result holds for any chordal graph $G'$ with at most $n-1$ vertices and any spanning subgraph $H'\subseteq G'$ satisfying $\Mad(H')\leq d$. 

If $G$ has a vertex $v \in V(G)$ such that $\deg_{G}(v) +2 \deg_{H}(v) < k$, let $G' = G-v$ and $H'=H-v$. Note that $G'$ is chordal, $\omega(G') \leq \omega$, and $H'$ is a spanning subgraph of $G'$ satisfying $\Mad(H')\leq \Mad(H)\leq d$. By induction hypothesis, let $\phi'$ be a 2-backbone $k$-colouring of $(G',H')$.
By the assumption on the degree of $v$ in $G$ and in $H$, at least one colour of $\{1,\dots,k\}$ remains valid to colour $v$ and extend $\phi'$ into a 2-backbone $k$-colouring of $(G,H)$. This shows that the result holds if $G$ contains such a vertex.
Consequently, let us assume that, for every vertex $v \in V(G)$, that
\begin{equation}
\label{eq:lowerbound_vertexdegree}
\deg_G(v) +2 \deg_H(v)  \geq k > (1+\epsilon)\omega+c_{\epsilon,d} -1. 
\end{equation}
In the sequel, we prove that this case cannot occur, by the choice of $c_{\epsilon,d}$, which concludes the proof.

For an ordering $(v_1,\ldots,v_n)$ over $V(G)$, define $X_i=\{v_j\mid 1\leq j \leq i\}$ and  $W_i=N_G[X_i]$. As shown in~\cite[Lemma~14]{BHP23}, there exists an ordering $(v_1,\ldots,v_n)$ over $V(G)$ such that, for all integers $1\leq i\leq n$:
\begin{gather}
    \deg_G(v_i) \leq i+\omega-2\text{, and} \label{eq:upperbound_vertexdegree}\\
    |W_i|\leq 2i + \omega-1. \label{eq:upperbound_unionofdegrees}
\end{gather} 
From~\eqref{eq:lowerbound_vertexdegree} and~\eqref{eq:upperbound_vertexdegree}, it follows that, for every integer $1\leq i\leq n$:
\begin{align*}
    2\sum\limits_{1 \leq j \leq i} \deg_H(v_j)  &=  \sum\limits_{1 \leq j \leq i} 2\deg_H(v_j)
    > \sum\limits_{1 \leq j \leq i} \left((1+\epsilon)\omega+c_{\epsilon,d} -1- \deg_G(v_j)\right)\nonumber\\
    &\geq  \sum\limits_{1 \leq j \leq i} \left((1+\epsilon)\omega+c_{\epsilon,d} - 1 - (j+\omega-2)\right)\\
    &= \sum\limits_{1 \leq j \leq i} \left(\epsilon \omega + c_{\epsilon,d} - j+1\right)\\
    &= i\epsilon \omega + ic_{\epsilon,d} - \frac{i(i-1)}{2}.
    \addtocounter{equation}{1}\tag{\theequation}
    \label{eq:lb_degvj}
\end{align*}

On the other hand, since $H$ satisfies $\Mad(H)\leq d$, we have that $2|E(H')|\leq d|V(H')|$, for any $H'\subseteq H$. Then, note that~\eqref{eq:upperbound_unionofdegrees} implies
\begin{equation}
    2\sum\limits_{1 \leq j \leq i} \deg_H(v_j) \leq 2\Big(|E(H[W_i])| + |E(H[X_i])|\Big)
    \leq d\Big(|W_i| + |X_i|\Big)
    \leq d(3i + \omega-1). \label{eq:ub_degvj}
\end{equation}

Hence, by~\eqref{eq:lb_degvj} and~\eqref{eq:ub_degvj}, for every integer $1\leq i\leq n$, we have $i\epsilon \omega + ic_{\epsilon,d} - \frac{i(i-1)}{2} < d(3i + \omega-1)$,
which implies that $c_{\epsilon,d}  < 3d - \epsilon \omega + \frac{d\omega}{i} + \frac{(i-1)}{2} -\frac{d}{i}$.
Recall that $n\geq \frac{d}{\epsilon}$. Therefore, applied $i = \left\lceil \frac{d}{\epsilon} \right\rceil$, the previous inequality implies that $c_{\epsilon,d}  < 3d+ \frac{d}{2\epsilon}$.
This contradicts our choice of $c_{\epsilon,d}$.
\end{proof}

As mentioned in the introduction, we thus obtain the following when $H$ is a forest.

\begin{corollary}
    \label{cor:upbound_chordal_trees}
    If $G$ is chordal and $H$ is a forest, then $\BBC_2(G,H)\leq \omega(G)+\Ocal(\sqrt{\omega(G)})$.
\end{corollary}
\begin{proof}
    It is a direct consequence of \Cref{thm:mad} as $\Mad(T)<2$ holds for any tree $T$.
\end{proof}

\subsection{Backbones without 4-cycles}

This section is devoted to the proof of Theorem~\ref{thm:3half}.

\thmcfourfree*

\begin{proof}
Assume first that $\omega(G)=\omega$ is odd.
There exists (see~\cite[Lemma 12]{BHP23}) a tree-decomposition $(T,{\cal X})$ of $G$ such that $T$ is rooted in $r\in  V(T)$ and, $|X_v|=\omega$ for all $v \in V(T)$ and, for every $v \in V(T) \setminus \{r\}$ with parent $p$, $|X_p\setminus X_v|=|X_v\setminus X_p|=1$.

In~\cite{BHP23} (proof of Theorem 7), it was further shown that there exists a proper colouring $\phi:V \rightarrow \{1,\cdots,\frac{\omega+3}{2}\}$ such that, for each $t \in V(T)$, every colour appears at most twice in $X_t$ and at most three colours appear uniquely in $X_t$. Moreover, using the fact that $H$ is $C_4$-free, two vertices that are adjacent in $H$ are assigned distinct colours.

This implies that $V(G)$ can be partionned into $k\leq \frac{\omega+3}{2}$ induced forests $F_1,\cdots,F_k$. For every $1 \leq i \leq k$, proper colour the vertices of $F_i$ using colours $3i-1$ and $3i-2$. This is clearly a proper colouring of $G$. Moreover, for every $\{u,v\} \in E(H)$, $u$ and $v$ are in distinct forests and so, their colours differ by at least $2$. 
Hence, $\BBC_2(G,H)\leq 3\frac{\omega+3}{2}-1=\frac{3\omega+7}{2}$.

Finally, if $\omega$ is even, we apply previous paragraph to $G-I$ (where $I$ is an independent set with $\omega(G-I)=\omega(G)-1$) and add two extra colours for $I$. We thus obtain 
\[
\BBC_2(G,H)\leq \frac{3(\omega-1)+7}{2}+2=\frac{3\omega+8}{2},
\]
as desired.
\end{proof}

\section{Further Research}
\label{sec:further}

In this paper, we prove several evidences for the following conjecture, which is still open.

\begin{conjecture}
    Let $G$ be a chordal graph and $H$ be a spanning forest of $G$. Then 
    \[
    \BBC_2(G,H)\leq \omega(G) + \Ocal(1).
    \]
\end{conjecture}

We actually believe that even the following much stronger one holds.

\begin{conjecture}
    There exists a function $f: \N^2 \to \N$ such that the following holds. For every chordal graph $G$ and every subgraph $H$ of $G$, if $\Mad(H) \leq d$ then 
    \[
    \BBC_q(G,H) \leq \omega(G) + f(q,d).
    \]
\end{conjecture}

When $H$ is bipartite, the following problem is open. The fact that $\frac{5}{3} \leq \gamma \leq 2$ follows from \eqref{eq:trivialbounds} and Proposition~\ref{prop:example_chordal_bip_53w}.

\begin{problem}
    Find the infimum of all values $\frac{5}{3} \leq \gamma \leq 2$ for which $X$ is infinite, where
    \[X=\{ \omega(G) \mid \exists G \text{ chordal},H \text{ bipartite}, H\subseteq G, \BBC_2(G,H) \geq \gamma \cdot \omega(G)\}.\]
\end{problem}

\bibliographystyle{abbrv}
\bibliography{BBC-Chordal-LAGOS2025}

\begin{thebibliography}{10}

\bibitem{AHK+07}
K.~I. Aardal, S.~P.~M. Hoesel, A.~M. C.~A. Koster, C.~Mannino, and A.~Sassano.
\newblock Models and solution techniques for frequency assignment problems.
\newblock {\em Annals of Operations Research}, 153(1):79--129, 2007.

\bibitem{appel1977}
K.~Appel and W.~Haken.
\newblock Every planar map is four colorable. {P}art {I}: {D}ischarging.
\newblock {\em Illinois J. Math.}, 21(3):429--490, 09 1977.

\bibitem{appel1977a}
K.~Appel, W.~Haken, and J.~Koch.
\newblock Every planar map is four colorable. {P}art {II}: {R}educibility.
\newblock {\em Illinois J. Math.}, 21(3):491--567, 09 1977.

\bibitem{AHS.15}
J.~Araujo, F.~Havet, and M.~Schmitt.
\newblock Steinberg-like theorems for backbone colouring.
\newblock {\em Electronic Notes in Discrete Mathematics}, 50:223--229, 2015.

\bibitem{AAC+22}
C.~S. Araújo, J.~Araújo, A.~Silva, and A.~A. Cezar.
\newblock Backbone coloring of graphs with galaxy backbones.
\newblock {\em Discrete Applied Mathematics}, 323:2--13, 2022.

\bibitem{ACT24}
J.~Araújo, R.~Castro, and A.~Talon.
\newblock About directed backbone colourings of graphs.
\newblock In {\em Annals of the 11th Latin American Workshop on Cliques in Graphs}, page~26, 2024.

\bibitem{BHP23}
S.~Bessy, F.~Havet, and L.~Picasarri-Arrieta.
\newblock Dichromatic number of chordal graphs.
\newblock {\em Graphs and Combinatorics}, 41(4):81, 2025.

\bibitem{BM.book}
J.~A. Bondy and U.~S.~R. Murty.
\newblock {\em Graph Theory}.
\newblock Springer, 2008.

\bibitem{broersma2009backbone}
H.~Broersma, B.~Marchal, D.~Paulusma, and A.~Salman.
\newblock Backbone colorings along stars and matchings in split graphs: their span is close to the chromatic number.
\newblock {\em Discussiones Mathematicae Graph Theory}, 29(1):143--162, 2009.

\bibitem{Broersma2003}
H.~J. Broersma, F.~V. Fomin, P.~A. Golovach, and G.~J. Woeginger.
\newblock Backbone colorings for networks.
\newblock In H.~L. Bodlaender, editor, {\em Graph-Theoretic Concepts in Computer Science}, pages 131--142. Springer, 2003.

\bibitem{BFGW07}
H.~J. Broersma, F.~V. Fomin, P.~A. Golovach, and G.~J. Woeginger.
\newblock Backbone colorings for graphs : tree and path backbone.
\newblock {\em Journal of Graph Theory}, 55(2):137--152, 2007.

\bibitem{BFY.03}
H.~J. {Broersma}, J.~{Fujisawa}, and K.~{Yoshimoto}.
\newblock Backbone colorings along perfect matchings.
\newblock Technical Report 1706, University of Twente, 2003.

\bibitem{CHSS13}
V.~Campos, F.~Havet, R.~Sampaio, and A.~Silva.
\newblock Backbone colouring: Tree backbones with small diameter in planar graphs.
\newblock {\em Theoretical Computer Science}, 487:50--64, 2013.

\bibitem{diestel2017}
R.~Diestel.
\newblock {\em Graph Theory}.
\newblock Graduate texts in mathematics. Springer, 5th edition, 2017.

\bibitem{HK.14}
F.~Havet and A.~D. King.
\newblock List circular backbone colouring.
\newblock {\em Discrete Mathematics and Theoretical Computer Science}, 16(1):89--104, 2014.

\bibitem{HKLT14}
F.~Havet, A.~D. King, M.~Liedloff, and I.~Todinca.
\newblock ({C}ircular) backbone colouring: Forest backbones in planar graphs.
\newblock {\em Discrete Applied Mathematics}, 169:119--134, 2014.

\bibitem{JT95}
T.~R. Jensen and B.~Toft.
\newblock {\em Graph Coloring Problems}.
\newblock Wiley-Interscience, 1995.

\bibitem{Karp72}
R.~Karp.
\newblock Reducibility among combinatorial problems.
\newblock In R.~Miller and J.~Thatcher, editors, {\em Complexity of Computer Computations}, pages 85--103. Plenum Press, 1972.

\bibitem{MR01}
M.~Molloy and B.~Reed.
\newblock {\em {Graph Colouring and the Probabilistic Method}}.
\newblock Springer, 2001.

\bibitem{salman2006lambda}
A.~N.~M. Salman.
\newblock $\lambda$-backbone coloring numbers of split graphs with tree backbones.
\newblock In {\em Proceeding of The Second IMT-GT 2006 Regional Conference on Mathematics, Statistics and Applications}, pages 43--47, 2006.

\bibitem{turowski2015}
K.~Turowski.
\newblock Optimal backbone coloring of split graphs with matching backbones.
\newblock {\em Discussiones Mathematicae Graph Theory}, 35(1):157--169, 2015.

\end{thebibliography}

\end{document}